\documentclass[12pt]{article}

\title{Functional limit theorems for \\ L\'evy processes 
 satisfying Cram\'er's condition} 
 
\author{M\'aty\'as Barczy\thanks{Department of Applied Mathematics and Probability, Faculty of Informatics, University of Debrecen, Pf.12, H-4010 Debrecen, Hungary. Email:
barczy.matyas@inf.unideb.hu}  \and
Jean Bertoin\thanks{Laboratoire de Probabilit\'es et Mod\`eles Al\'eatoires, 
UPMC, 4 place Jussieu, F-75252 Paris cedex 05, France. Email:
jean.bertoin@upmc.fr}  }

\date{}

\usepackage{amsfonts,amsmath,amssymb,bbm}
\usepackage{setspace} 
\def\proof{\noindent{\bf Proof:}\hskip10pt}        
\def\QED{\hfill $\Box$}

\font\tenmath=msbm10 scaled 1200
\font\sevenmath=msbm7 scaled 1200
\font\Phiivemath=msbm5 scaled 1200
\newfam\mathfam \textfont\mathfam=\tenmath
\scriptfont\mathfam=\sevenmath \scriptscriptfont\mathfam=\Phiivemath

\def \\ { \cr }
\def\R{\mathbb{R}}
\def \1{1 \mkern -6mu 1} 

\def\E{\mathbb{E}}
\def\P{\mathbb{P}}

\def\D{\mathbb{D}}

\def \f{{\mathcal F}}

\def \up{^{\uparrow}}
\def \da{^{\downarrow}}
\def \e{{\rm e}}
\def \d{{\rm d}}

\def \p{{\mathcal P}}
\def \q{{\mathcal Q}}

\newtheorem{theorem}{Theorem}

\newtheorem{lemma}{Lemma}
\newtheorem{corollary}{Corollary}
\begin{document}

\maketitle

\begin{abstract}
We consider  a L\'evy process that  starts from
$x<0$ and conditioned on having a positive maximum. When Cram\'er's condition holds, we provide two weak limit theorems as $x\to -\infty$ for the law of  the (two-sided) path shifted at the first instant when it enters $(0,\infty)$, respectively shifted at the instant when its overall maximum is reached. The comparison of these two asymptotic results yields some interesting identities related to time-reversal, insurance risk, and self-similar Markov processes. \newline

\textbf{2000 Mathematics Subject Classification: 60G51, 60G18, 60B10.} 

{\bf Key words:} L\'evy process, Cram\'er's condition, self-similar Markov process.

 \end{abstract}

\begin{section}{Introduction}

This work should be viewed as a continuation of \cite{BS} which dealt with certain weak limit theorems  for  real-valued L\'evy processes having integrable ladder heights. Recall that this includes the case when the L\'evy process has a finite and positive first moment and rules out the case when the L\'evy process drifts to $-\infty$. The main result of \cite{BS} is that the distribution of such a L\'evy process started from $x$ and shifted at the time of its first entrance in $[0,\infty)$, has a weak limit as $x\to -\infty$. The latter  is the law  of a process indexed by the entire real line, which remains negative for negative times, crosses $0$ at the instant $0$ and then evolves as the original L\'evy process.  This two-sided process fulfills a remarkable property of spatial stationarity and yields in particular a construction of self-similar Markov processes in $(0,\infty)$ starting from the entrance boundary point $0$  via a Lamperti-type transformation, see Corollaries 3 and 4 in \cite{BS}.

In the present paper, we consider a L\'evy process $\xi$ which fulfills Cram\'er's condition, namely we suppose 
$$\hbox{there exists $\theta>0$ such that }\ E(\exp(\theta \xi_1))=1\,.$$
Recall that {\it a fortiori} $\xi$ drifts to $-\infty$ and hence the ladder heights are now defective. 
It is well-known that if we introduce the exponentially tilted probability measure (also known as Esscher transform)
$$\tilde P(\Lambda)= E(\exp(\theta \xi_t) {\bf 1}_{\Lambda})\,,\qquad \Lambda\in\f_t\,,$$
where $(\f_t)$ stands for the natural filtration of $\xi$, then $\xi$ remains a L\'evy process under $\tilde P$ and has a finite and positive first moment whenever
$$E(|\xi_1|\exp(\theta \xi_1))<\infty\,.$$

Applying classical techniques of changes of probability which have been developed by Asmussen \cite{As1} in the setting of random walks (see also \cite{As2}  and  \cite{Ky} for  L\'evy processes),  it is then easy to deduce from \cite{BS} that the initial L\'evy process started from $x$, {\it conditioned to visit} $[0,\infty)$ and shifted in time to its first entrance in $[0,\infty)$, has a weak limit as $x\to -\infty$ which can be expressed in terms of the spatially stationary version of the exponentially tilted L\'evy process.

The starting point of this work lies in the observation that one can also investigate the asymptotics of the preceding conditional laws using a path decomposition at the instant when the process reaches its maximum, instead of the first entrance time in $[0,\infty)$. We thus have two equivalent descriptions (see Theorems \ref{T1} and \ref{T2} below) of the same two-sided process  and the comparison  yields interesting identities in distribution.
One of the main applications concerns again self-similar Markov processes and more precisely  their recurrent extensions in $[0,\infty)$. We will derive a path-decomposition {\it \`a la} Williams \cite{DW} (i.e. at the instant when the maximum is attained) under It\=o's measure of their excursions away from $0$. The latter involves through a Lamperti's type transformation two independent versions of the underlying L\'evy process conditioned to stay positive and to stay negative, respectively. 
In a different direction, we also derive an application to risk theory. Namely, we determine the asymptotic distribution of the total time when an insurance is indebted for large initial reserve and conditionally on the rare event that ruin occurs. 

Very recently, Griffin and Maller \cite{Gr, GM} have investigated a closely related topic. Basically, 
they deal with L\'evy processes whose L\'evy measure belongs to the convolution equivalent class 
${\mathcal S}^{(\alpha)}$ and also consider asymptotic path-decompositions conditionally on the rare event that the process crosses a large level. There are however several  fundamental differences, besides the simple fact that the classes of  L\'evy processes which are considered in \cite{Gr, GM} and here are not the same. First, our results involve a single limit (we let the starting point of the L\'evy process tend to $-\infty$) whereas Griffin and Maller have a two-step procedure (they consider the first passage time above $u-x$ given that the process exceeds level $u$, and let first $u\to \infty$ and then $x\to \infty$; cf. Theorem 3.1 in \cite{GM} or Theorem 6.2  in \cite{Gr}). 
A second key difference is that Griffin and Maller describe processes forward in time, while our description is two-sided and uses  time reversal for one part. 

The plan of the rest of this paper is as follows. Section 2 is devoted to notations and preliminaries on conditioning and Cram\'er's hypothesis. The two main limit theorems are established in Section 3, 
and Section 4 is dedicated to some applications. 

 \end{section}

 \begin{section}{Preliminaries}
 
  We introduce notations that may look at first sight more complicated than needed, but will turn out to be useful  for our purposes.

\subsection{Canonical notation}
We work with generic c\`adl\`ag paths, possibly with finite lifetime, which can be either one-sided
$\omega: [0,\infty)\to \R\cup\{-\infty\}$ or two-sided $\bar \omega: \R\to \R\cup\{-\infty\}$, where $-\infty$ is viewed as an isolated extra point added to $\R$. We shall frequently view a two-sided path $\bar \omega$ as a pair of one-sided paths $(\bar \omega(t), t\geq 0)$ and $(-\bar \omega(-t-), t\geq 0)$ put back-to-back.
We write respectively $\Omega$ and $\bar \Omega$ for the spaces of such paths
 and endow both with the  Skorokhod-Lindvall topology. 

It will sometimes be convenient to extend a one-sided path $\omega\in \Omega$ to $\R$ by declaring that $\omega(t)=-\infty$ for all $t<0$. More generally, by a slight abuse of notation, we shall often consider a real valued path $\omega': [a,b)\to \R$ defined on some interval  $[a,b)$ with $-\infty<a < b \leq \infty$ as an element of $\Omega$ (provided that $a\geq 0$) or of $\bar \Omega$ by agreeing  that $\omega'(s)=-\infty$
for $s\not\in [a,b)$. We then refer to $[a,b)$ as  the life-interval of $\omega'$.
In this direction, we will use the notation $\sup \omega'=\sup_{a\leq t < b} \omega'(t)$
and $\inf \omega'=\inf_{a\leq t < b} \omega'(t)$ respectively for the supremum and the infimum over the life-interval of $\omega'$.

  We  denote by 
 $\xi=(\xi_t)$ the canonical process on $\Omega$ or $\bar \Omega$, viz. $\xi_t(\omega)=\omega(t)$,
 and by $(\f_t)$ the canonical filtration, that is $\f_t$ is the smallest sigma-field on $\Omega$
 or $\bar \Omega$
 such that $\xi_s$ is measurable for all $s\leq t$.
   We consider a probability measure $P$ on $\Omega$ under which $\xi$ is a L\'evy process.
 By this, we mean that under $P$,   the life-interval is $[0,\infty)$ a.s.  (i.e. $\xi_t\in\R$ for all $t\geq 0$, $P$-a.s.), and for every $s,t\geq 0$, 
the increment $\xi_{t+s}-\xi_t$ is independent of $\f_t$ and has the same law as $\xi_{s}$.
For every $x\in\R$, we write $P_x$ for the law of $x+\xi$ under $P$; in other words $P_x$ is the law of the L\'evy process started from $x$ at time $0$.

\subsection{Conditioning to stay negative}

We shall assume throughout  this work that $0$ is regular for both half-lines. This means that if we introduce the entrance times
$$T=\inf\{t>0: \xi_t<0\} \quad \hbox{and} \quad  \tau=\inf\{t>0: \xi_t>0\}\,,$$
then
\begin{equation}\label{Eq1}
T=\tau=0\,,\qquad P\hbox{-a.s.}
\end{equation}
The hypothesis \eqref{Eq1} is merely taken for granted for the sake of convenience, and it should not be difficult for an interested reader to adapt the arguments of this work and deal with the situation where $0$ is irregular for $(-\infty, 0)$ or for $(0,\infty)$.
Note that \eqref{Eq1} in particular rules out the compound Poisson case.
Further essential assumptions will be made in due time.

We  also suppose that the L\'evy process drifts to $-\infty$, in the sense that
$$\lim_{t\to \infty}\xi_t=-\infty\qquad P\hbox{-a.s.}$$
Since for $x<0$, there is a positive probability that the L\'evy process started from $x$ always remains negative, we can introduce the conditional  law 
$$P\da_x:=P_x(\cdot \mid \sup \xi <0)\,.$$
According to Theorem 2 of Chaumont and Doney \cite{CD}, \eqref{Eq1} ensures that $P\da_x$ has a weak limit as $x\to 0-$ which we denote by
$P\da$, or sometimes also by $P\da_0$ when it is more convenient.  
This is the law of a Markov process started from $0$ and which is negative at positive times; its semigroup is that of the L\'evy process conditioned to stay negative.

\subsection{Cram\'er's condition}
We also assume that Cram\'er's condition is fulfilled, that is 
\begin{equation}\label{Eq2}
\exists \theta >0 \quad\hbox{ such that }\quad E(\exp(\theta \xi_t))=1
 \end{equation}
and furthermore  that
\begin{equation}\label{Eq3}
 E(|\xi_t|\exp(\theta \xi_t))<\infty\,,
\end{equation}
for some (and then all) $t>0$.
Recall that this automatically implies that the L\'evy process drifts to $-\infty$.
Cram\'er's condition \eqref{Eq2} entails that the exponential L\'evy process $\exp(\theta \xi_t)$ is a $P$-martingale, and this enables us to introduce the tilted probability measure
$$\tilde P(\Lambda)= E(\exp(\theta \xi_t) {\bf 1}_{\Lambda})\,,\qquad \Lambda\in\f_t\,,$$
where $t\geq 0$ is arbitrary. It is well-known that under $\tilde P$, the canonical process $\xi$ is still a L\'evy process which now drifts to $\infty$, i.e. $\tilde P(\lim_{t\to \infty} \xi_t=\infty)=1$; 
one should think of $\tilde P$ as the law of the L\'evy process conditioned to drift to $\infty$.

We now state a slight extension of Cram\'er's classical estimate (cf. \cite{BD}) which will be useful later on.
Denote by 
$$\sigma=\inf\{t\geq 0: \xi_t \vee \xi_{t-}=\sup \xi \}$$
the (first)  instant when the canonical process reaches its overall supremum. It is well-known that when \eqref{Eq1} holds, then  $\xi$
hits  its maximum continuously and precisely at time $\sigma$, i.e.
$$P(\xi_{\sigma}=\xi_{\sigma-}=\sup \xi \hbox{ and }\xi_t\vee \xi_{t-}<\sup \xi \hbox{ for all }t\neq \sigma)=1\,;$$
see for instance Corollary 3.11 in \cite{Ky}.

\begin{lemma}\label{L0} Assume that  \eqref{Eq1},  \eqref{Eq2} and \eqref{Eq3} hold. There is a constant $C>0$ such that for every $t\geq0$, we have
$$\lim_{x\to -\infty} \e^{-\theta x}P_x(\sup \xi  \geq 0, \sigma\geq t)=C\,.$$
\end{lemma}
We stress that the constant $C$ in Lemma \ref{L0} can be represented in several different ways, see \cite{BD}.

\proof For $t=0$, this is just Cram\'er's estimate, so assume that $t>0$. 
Recall that $\tau$ denotes the first entrance time in $(0,\infty)$, so for all $x\in\R$
$$\e^{-\theta x}P_x(\sup \xi  \geq 0, \sigma< t) \leq \e^{-\theta x} P_x(\tau<t)\,.$$
We then use the exponential tilting  and the optional sampling theorem
to express the right-hand side as
$$\tilde E_x(\e^{-\theta \xi_{\tau}}, \tau<t) \leq \tilde P_x(\tau<t)\,,$$
where $\tilde P_x$ denotes the law of $x+\xi$ under $\tilde P$.
Plainly the preceding quantity tends to $0$ as $x\to -\infty$, which yields our claim. \QED

Next, just as in the preceding section, we can introduce for $x>0$ the conditional law 
$$\tilde P\up_x= \tilde P_x(\cdot \mid \inf \xi  >0)$$
 together with its weak limit 
$$\tilde P\up =\tilde P\up_0= \lim_{x\downarrow 0} \tilde P\up_x\,.$$

Recall also that $\tau$ denotes the first instant when $\xi$ enters $(0,\infty)$ and that due to \eqref{Eq3},
 under the tilted law $\tilde P$, the L\'evy process  $\xi$ has a finite positive expectation (and hence its ladder heights are integrable).
 In this situation,   the distribution  under $\tilde P_z$ of the  pair  $(-\xi_{\tau-}, \xi_{\tau})$ formed by the under-shoot and over-shoot, has a weak limit as $z\to -\infty$, 
\begin{equation}\label{Eq4}
\lim_{z\to -\infty} \tilde P_z(-\xi_{\tau-} \in \d x, \xi_{\tau}\in \d y) = \tilde \rho(\d x , \d y)\,, \qquad
x,y \in\R_+\,,
\end{equation}
 where  the limiting probability measure $\tilde \rho$ is described by Equation (5) in \cite{BS}; see
 for instance Lemma 3 there.

 It has been pointed out in  \cite{BS} that the weak convergence \eqref{Eq4} of the under-shoot and over-shoot   can be extended to paths.
  In this direction,
 we first introduce a probability measure $\tilde\p$ on
       $\bar\Omega$ by
      $$
               \tilde\p(\d \bar\omega) = \tilde\p(\d \omega, \d \omega')
           = \int_{\R_+\times\R_+} \tilde\rho(\d x, \d y)
             \tilde P_x^{\uparrow}(\d\omega) \tilde P_y(\d\omega'),
       $$
       where a two-sided path $\bar \omega$ is identified with a pair of one-sided
       paths $(\omega,\omega')$ as it has been explained in the preliminaries.
     In other words,  under $\tilde\p$ the distribution of $(-\xi_{0-},\xi_0)$ is $\tilde\rho$,
       and conditionally on $(-\xi_{0-},\xi_0)=(x,y)$, the processes $(-\xi_{-t-},t\geq 0)$
       and $(\xi_t,t\geq 0)$ are independent with laws $\tilde P_x^{\uparrow}$ and
       $\tilde P_y$, respectively.
       According to Theorem 3 in \cite{BS}, we have 
\begin{equation}\label{Eq5}
{\mathcal L}(\xi_{\tau +\cdot},\tilde P_x)\quad \Longrightarrow \quad {\mathcal L}(\xi,\tilde \p)
\qquad \hbox{as }x\to-\infty\,,
\end{equation}
where the notation ${\mathcal L}(\eta,Q)$ refers to the law of the process $\eta$ under a probability measure $Q$ and $\Rightarrow$ to weak convergence  in the sense of Skorokhod-Lindvall, and 
by a slight abuse of notation, we view the process shifted in time at its first entrance in $(0,\infty)$,
$\xi_{\tau + \cdot} : t\mapsto \xi_{\tau +t}$,  as a two-sided process  by agreeing that $ \xi_{\tau +t}=-\infty$ for $t<-\tau$. For being precise, we note that Theorem 3 in \cite{BS} states
      that the weak convergence \eqref{Eq5}  holds on $\D([b,\infty))$ for any fixed $b\in\R$.
      However, it is equivalent that this holds on $\D(\R)$, see for instance \cite{WW}.

\end{section}

 \begin{section}{Limits under conditioning to have a large maximum}
 Our main goal is to investigate the asymptotic behavior of the initial L\'evy process in the neighborhood of the instant when it reaches its maximum, conditionally on the event that its maximum is abnormally large. We shall first tackle this question via exponential tilting, relying on \eqref{Eq5}, and then independently using classical techniques of path decomposition at the instant when the maximum is reached. 
 
   \subsection{Shifting at the first entrance time in $[0,\infty)$}
   We start by recalling that for every $x\in\R$, the process
   $(\exp(-\theta \xi_t), t\geq 0)$ is a $\tilde P_x$-martingale, and that
   for every $(\f_t)$-stopping time $S$ which is  finite $\tilde P_x$-a.s., there is the identity
   $$P_x(\Lambda, S<\infty)= \e^{\theta x} \tilde E_x(\exp(-\theta \xi_S) {\bf 1}_{\Lambda})\,,\qquad \Lambda\in\f_S\,.$$
   This enables us to introduce a probability measure $\p$ on $\bar \Omega$ by converse exponential tilting of the limit law $\tilde \p$. Specifically we set
  \begin{equation} \label{Eq6}  
  \p(\Lambda, S<\infty)= c(\theta) \tilde {\mathcal E}(\exp(-\theta \xi_S) {\bf 1}_{\Lambda})\,,\qquad \Lambda\in\f_S\,,
  \end{equation}
    where $S\geq 0$ is an $(\f_t)$-stopping time which is $\tilde{\p}$-a.s. finite, $\tilde {\mathcal E}$ refers to the expectation under $\tilde \p$, and
    $c(\theta)>0$ is the normalizing constant, i.e. 
    $$\frac{1}{c(\theta)}=\int_{\R_+\times \R_+} \e^{-\theta y}Ê\tilde \rho(\d x, \d y)=\tilde{\mathcal E}(\exp(-\theta \xi_0))\,.$$

Recalling the construction of $\tilde \p$ in the preceding section, we obtain an equivalent description of $\p$. First,
\begin{equation}\label{Eq7}
\p(-\xi_{0-}\in \d x, \xi_0\in \d y)  = c(\theta) \e^{-\theta y}Ê\tilde \rho(\d x, \d y):= \rho(\d x, \d y) \,,\qquad x,y \in\R_+\,,
\end{equation}
and  then, under the conditional law
 $\p(\cdot \mid -\xi_{0-}=x, \xi_0=y)$, the processes  $(-\xi_{-t-}, t\geq 0)$ and $(\xi_t, t\geq 0)$ are independent with respective laws
 $\tilde P_x\up$ and $P_y$.  Note that $\lim_{t\to\pm\infty}\xi_t=-\infty$, $\p$-a.s

  Recall from Corollary 3 in \cite{BS} that $\tilde{\p}$ fulfills the spatial stationarity property: for all $y\in\R$
 \begin{equation}\label{Eq8}
 {\mathcal L}(\xi_{\tau_y +\cdot},\tilde {\p})=
 {\mathcal L}(y+\xi,\tilde{\p})\,,
 \end{equation}
  where $\tau_y=\inf\{t: \xi_t>y\}$ denotes the first passage time above $y$.
This yields the spatial quasi-stationarity property for $\p$, namely
for  every $y>0$, there is the identity
\begin{equation}\label{Eq9}
{\mathcal L}(\xi_{\tau_y +\cdot},\p(\cdot \mid \sup \xi >y))=
 {\mathcal L}(y+\xi,\p)\,.
 \end{equation}
 Indeed,  specifying  \eqref{Eq6}  for $S=\tau_y$, we find in particular using \eqref{Eq8} that
$\p(\sup \xi >y)=\e^{-\theta y}$, and then \eqref{Eq9} follows from \eqref{Eq8} and \eqref{Eq6} applied to $S=\tau_y+t'$ for $t'\geq 0$ arbitrarily large.
 
Similarly, the following limit theorem derives readily from Theorem 3 in \cite{BS} .
   
   \begin{theorem}\label{T1} Assume \eqref{Eq1},  \eqref{Eq2} and \eqref{Eq3}.
 We have  that
 $${\mathcal L}(\xi_{\tau +\cdot},P_x(\cdot \mid \sup \xi >0))\quad \Longrightarrow \quad \p \qquad \hbox{as }x\to-\infty\,.
$$
\end{theorem}

\proof Fix $a<0<b$, consider a bounded continuous functional $\Phi: \bar \Omega \to \R$ which only depends on the restriction of $\bar \omega$ to the time-interval $[a,b]$, and write $\eta=\Phi((\xi_{\tau +t})_{a\leq t \leq b})$. We have for all $x\in \R$ that
$$E_x(\eta \mid \sup \xi >0)
= \frac{E_x( \eta, \tau <\infty)}
{P_x( \sup \xi >0)}
=\frac{\e^{\theta x}}
{P_x( \sup \xi >0 )}\, \tilde E_x( \eta  \exp(-\theta \xi_{\tau+b}))\,.
$$
On the one hand,  \eqref{Eq5} shows that
$$\lim_{x\to -\infty} \tilde E_x( \eta  \exp(-\theta \xi_{\tau+b}))
=\tilde {\mathcal E} \left( \Phi((\xi_{t})_{a\leq t \leq b})  \exp(-\theta \xi_{b})\right)\,.$$
On the other hand, we know from Cram\'er's estimate  that
$$P_x( \sup \xi >0 ) \sim C \e^{\theta x}\qquad \hbox{ as }x\to -\infty\,,$$
where $C$ is the expected value of $\exp(-\theta \tilde \gamma)$ with a variable $\tilde \gamma$ distributed as the stationary overshoot for the tilted L\'evy process (i.e. as the second marginal of $\tilde \rho$);   see Remark (2) in \cite{BD}. Hence we have $c(\theta)C=1$, which completes the proof. \QED

\subsection{Shifting at the time of the maximum}
We start by introducing another probability measure $\q$ on $\bar{\Omega}$
as follows. Consider two independent processes $(\eta\da_t)_{t\geq 0}$ and 
$(\tilde \eta\up_t)_{t\geq 0}$ with respective laws $P\da$ and  $\tilde P\up$, together with an independent exponential variable $\epsilon$ with parameter $\theta$. 
Then define the two-sided process $(\eta_t)_{t\in\R}$ by
$$\eta_t=\left\{ 
\begin{matrix} \epsilon + \eta\da_t & \ \hbox{if }&t\geq 0\,, \\
\epsilon - \tilde\eta\up_{-t-}  & \ \hbox{if }&t\leq 0\,.\\
\end{matrix} \right.
$$
We denote the law of $\eta$ on $\bar \Omega$ by $\q$.
Observe that $\eta$ reaches its overall maximum at time $0$, and that this maximum is exponentially distributed, so that by the lack of memory of the exponential law, $\q$ also fulfills a 
spatial quasi-stationarity property. More precisely,
for  every $y>0$, there is the identity
$$
{\mathcal L}(\xi,\q(\cdot \mid \sup \xi >y))=
 {\mathcal L}(y+\xi,\q)\,.
$$
Note that this  is slightly different from \eqref{Eq9} as here there is no shift in time.

Recall that $\sigma=\inf\{t\geq 0: \xi_t = \sup  \xi \}$ denotes the  instant when the path reaches its overall maximum; we now claim a second weak limit theorem for the L\'evy process conditioned to have a large maximum.

\begin{theorem}\label{T2} Assume \eqref{Eq1},  \eqref{Eq2} and \eqref{Eq3}. We have that
 $${\mathcal L}(\xi_{\sigma +\cdot},P_x(\cdot \mid \sup \xi >0))\quad \Longrightarrow \quad \q
 \qquad \hbox{as }x\to-\infty\,.$$
\end{theorem}
It is well-known that under $P$, the 
post-supremum process $(\xi_{\sigma+t}-\sup \xi )_{t\geq 0}$ has the law $P\da$
and is independent of the pre-supremum process $(\xi_t)_{0\leq t<\sigma}$; see for instance Theorems 3.1 and 3.4 in
\cite{Be1}.
This reduces the proof of Theorem \ref{T2} to the following weaker statement.

\begin{lemma}\label{L1} Assume \eqref{Eq1},  \eqref{Eq2} and \eqref{Eq3}. The distribution of the process
$(\xi_{(\sigma-t)-})_{0\leq t<\sigma}$
under the conditional law
$P_x(\cdot \mid \sup \xi >0)$ 
converges weakly as $x\to -\infty$ to that of $(\epsilon-\tilde \eta\up_t)_{t\geq 0}$, where $\tilde \eta\up$ has the law $\tilde P\up$
and $\epsilon$ is an independent exponential variable with parameter $\theta$.
\end{lemma}

The rest of this section is devoted to the proof of Lemma \ref{L1}. In this direction, it is convenient to
denote by $\overleftarrow{\xi}$ the reversed pre-maximum path which has life-interval $[0,\sigma)$
 and is given by  
 $$\overleftarrow{\xi}_t= \sup \xi - \xi_{(\sigma-t)-}\,,Ê\qquad 0\leq t <\sigma\,.$$
To start with, we recall that the L\'evy process reflected at its running infimum $(\xi_t- \inf_{0\leq s\leq t}Ê\xi_s)_{t\geq 0}$ is Markovian under $P$. We denote by $n$ the It\=o measure of its excursions away from $0$, which is viewed as a sigma-finite measure on paths  killed when entering the negative half-line (i.e. at time $T$). Similarly, we write $\tilde n$ for the It\=o's excursion measure under the tilted law $\tilde P$. The connexion between the two excursion measures and the law $\tilde P\up$ is
that for every $t>0$ and every $\f_t$-measurable functional $F: \Omega\to \R_+$ we have
\begin{equation}\label{Eq10}
n(F(\xi) \e^{\theta \xi_t}, t<T)=\tilde n(F(\xi), t<T) = \tilde  E\up(F(\xi)/\tilde h(\xi_t))\,,
\end{equation}
where $\tilde h(x)=\tilde P_x(\inf \xi \geq 0)$. 
Indeed this follows readily from Corollary 1 in \cite{CD} and Theorem 3 in \cite{Ch}. 
More precisely,  the excursion measures $n$ and $\tilde n$ are only defined up to a constant factor that depends on the normalization of the local time at $0$, and hence we implicitly assume that  the normalization has been chosen such that \eqref{Eq10} holds.

We also denote  by $\kappa$
the Laplace exponent of the inverse local time at $0$ of the reflected L\'evy process (also called  the descending ladder times) which is given by 
\begin{equation}Ê\label{Eq11}
\kappa(q)=n(1-\e^{-qT})\,,\qquad q\geq 0\,,
\end{equation}
since it is well-known that the regularity assumption
\eqref{Eq1} implies that the descending ladder time process has zero drift.
Last, we introduce $k_t$, the killing operator at time $t$ on $\Omega$, so that the path
$\xi\circ k_t$ coincides with $\xi$ on the time-interval $[0,t)$ and takes the value $-\infty$
on $[t,\infty)$.

We are now able to express the law of the reversed pre-maximum process in terms of $n$ 
as follows.

\begin{lemma} \label{L2} Assume \eqref{Eq1}. Let $F$ be some nonnegative functional on $\Omega$. We have
$$E(F(\overleftarrow{\xi}))= \frac{1}{\kappa'(0)}\int_0^{\infty}n\left(F(\xi\circ k_t), t<T \right)  \d t \,.$$

\end{lemma}

\proof For $r>0$, denote by $P^r$ the law of the initial L\'evy process killed at rate $r$, i.e. the law of 
$\xi\circ k_{\epsilon(r)}$ under $P$, where $\epsilon(r)$ stands for an independent exponential variable with parameter $r$. It is a classical consequence of the duality lemma (Lemma II.2 in \cite{Be})  that under $P^r$, $\overleftarrow{\xi}$ has the same law as the post-minimum process $\underrightarrow{\xi}$, where
$\underrightarrow{\xi}_t = \xi_{\underline{\sigma} +t}-\min_{s<\epsilon(r)}\xi_s$ for $0\leq t < \epsilon(r)-\underline{\sigma}$ and $\underline{\sigma}$ is the a.s. unique instant at which the path reaches its overall minimum. 

Imagine that we mark times in $[0,\infty)$ using an independent Poisson process with intensity $r$.
Recall that since $0$ is regular for $(0,\infty)$, the set of times $t$ at which 
$\xi_t= \inf_{0\leq s\leq t}Ê\xi_s$ has zero Lebesgue measure, and hence contains no mark $P$-a.s. 
In this setting the post-minimum process $\underrightarrow{\xi}$ under $P^r$ can be identified as the
first excursion of the L\'evy process reflected at its running infimum away from $0$ which contains a mark and killed at this mark. It should then be plain from excursion theory and elementary properties of marked Poisson point processes that the law of $\underrightarrow{\xi}$ under $P^r$ coincides with
the distribution of the excursion conditioned to have at least one mark on its life-interval, and killed at this first mark. In other words, 
there is the identity
$$E^r(F(\underrightarrow{\xi}))= \frac{r}{\kappa(r)}\int_0^{\infty}n\left(F(\xi\circ k_t), t<T\right) \e^{-rt} \d t\,,$$
where $\kappa(r)=n(1-\e^{-rT})$ is the It\=o measure of the excursions containing at least one mark
according to \eqref{Eq11}. 
Our statement thus follows letting $r\to 0+$; recall that $\kappa'(0)<\infty$ since $\xi$ drifts to $-\infty$ under $P$. \QED

Next we introduce the potential operator for the L\'evy process killed when entering the negative half-line
$$Vg(x)=E_x\left(\int_0^{T}g(\xi_t)\d t\right)\,,\qquad x\geq 0\,,$$
where $g:\R \to\R$ denotes a generic measurable function. 

We now state a consequence of the classical Cram\'er's estimate for L\'evy processes.
\eject

\begin{lemma}\label{L3} Assume that  \eqref{Eq2} and \eqref{Eq3} hold. Let $g: \R\to\R$ be continuous with compact support
and set 
$g_x(y)=g(x+y)$ for every $x,y\in\R$. Then
$$\sup_{x,y\in\R}\e^{-\theta(x+y)}|Vg_x(y)|<\infty$$
and 
$$\lim_{x\to -\infty} \e^{-\theta x} Vg_x(y) = \frac{1}{\tilde m} \e^{\theta y} \tilde h(y) \int_{-\infty}^{\infty} \e^{-\theta z} g(z) \d z\,, \qquad y\in\R\,,$$
where $\tilde m=\tilde E(\xi_1)= E(\xi_1\exp(\theta \xi_1))$ and $\tilde h(y) = \tilde  P_y(\inf \xi >0)$.
\end{lemma}

\proof  An application of the strong Markov property at the first-passage time $T$ yields
$$Vg_x(y)= Ug_x(y) - E_y(Ug_x(\xi_{T}))\,,$$
where $U$ denotes the potential operator of the original L\'evy process, i.e.
$$Uf(y)=E\left(\int_0^{\infty} f(\xi_t+y)\d t\right).$$
Because $Ug_x(y)=Ug(x+y)$, we may rewrite the preceding as
$$Vg_x(y)= Ug(x+y) - E_y(Ug(x+\xi_{T}))\,.$$

Next denote by $\tilde U$ the potential operator of the  tilted L\'evy process,
$$\tilde Uf(y)=\tilde E\left(\int_0^{\infty} f(\xi_t+y)\d t\right) = E\left(\int_0^{\infty} \e^{\theta \xi_t}f(\xi_t+y)\d t\right)  \,,$$
and observe the identity
$$Uf(y) = \e^{\theta y} \tilde U \tilde f (y)\,, \quad \hbox{where } \tilde f(z)=\e^{-\theta z}f(z)\,.$$
As a consequence there is the upper-bound
$$\sup_{x,y\in\R}\e^{-\theta(x+y)}|Vg_x(y)|\leq \sup_{x,y\in\R}\e^{-\theta(x+y)}U|g|(x+y)\leq \sup_{z\in\R} \tilde U|\tilde g|(z)\,,$$
where $\tilde g$ is a continuous function with compact support. 
It follows readily that the right-hand side above is finite.

Recall also from the renewal theorem (see e.g. Theorem I.21 in \cite{Be}) that for every function $f$ which is continuous with compact support, we have
$$\lim_{x\to -\infty} \tilde U f(x+y)= \frac{1}{\tilde m} \int_{-\infty}^{\infty} f(z)\d z\,,\qquad y\in\R\,,$$
with $\tilde m=\tilde E(\xi_1)\in(0,\infty)$. As a consequence, we get that for every fixed $y\in\R$,
$$\lim_{x\to -\infty} \e^{-\theta x}Ug(x+y)= \frac{\e^{\theta y}}{\tilde m} \int_{-\infty}^{\infty} \e^{-\theta z} g(z)\d z\,.$$

Putting the pieces together, we finally conclude using dominated convergence that
$$\lim_{x\to -\infty} \e^{-\theta x} Vg_x(y) 
=  \frac{1}{\tilde m} \left(\int_{-\infty}^{\infty} \e^{-\theta z} g(z)\d z\right)
\left( \e^{\theta y}- E_y(\e^{\theta \xi_{T}})\right)\,.$$
Note that the last factor in the right-hand side can also be expressed as
$$\e^{\theta y} (1-\tilde P_y(T<\infty))=\e^{\theta y} \tilde P_y(T=\infty)\,,$$ 
since the stopping time $T$ is finite $P_y$-a.s.
 \QED

We are now able to tackle the proof of Lemma \ref{L1}.

\noindent {\bf Proof of Lemma \ref{L1}:}\hskip10pt  
Fix $t > 0$, consider a bounded continuous functional $F: \D[0,t]\to \R_+$
and a continuous function $g: \R_+\to \R_+$ with compact support. It will be convenient to agree that $g(y)=0$ for $y<0$.

We aim at estimating
$$
E_x\left(F(\overleftarrow{\xi}_s, 0\leq s \leq t) g(\sup \xi ), \sigma>t\right) 
= E\left(F(\overleftarrow{\xi}_s, 0\leq s \leq t) g(x+\sup \xi ), \sigma>t\right)
$$
when $x\to -\infty$.
We first use Lemma \ref{L2} and the notation $g_x(y)=g(x+y)$ to express this quantity as
\begin{eqnarray*}
& &\frac{1}{\kappa'(0)} \int_t^{\infty} n\left(F(\xi\circ k_t) g_x(\xi_s), s<T \right) \d s \\
&=& 
\frac{1}{\kappa'(0)} n\left(F(\xi\circ k_t) E_{\xi_t}\left(\int_0^{T} g_x(\xi_s)  \d s\right), t<T  \right)\\
&= &\frac{1}{\kappa'(0)} n\left(F(\xi\circ k_t) Vg_x({\xi_t}) , t<T \right)\,,
\end{eqnarray*}
where the first equality follows from the Markov property of the excursion measure $n$, and in the second, 
$V$ denotes the potential operator of the L\'evy process killed when it enters the negative half-line. 

Next we condition on $\sup \xi >0$ and apply \eqref{Eq10} to get 
\begin{eqnarray*}
& &E_x\left(F(\overleftarrow{\xi}_s, 0\leq s \leq t) g(\sup \xi ), \sigma>t\mid \sup \xi >0 \right) \\
&=& \frac{1}{\kappa'(0)P_x(\sup \xi >0)} n\left(F(\xi\circ k_t) Vg_x({\xi_t}) , t<T \right)\\
&=& \frac{\e^{\theta x}}{\kappa'(0)P_x(\sup \xi >0)} \tilde n\left(F(\xi\circ k_t) \e^{-\theta (x+\xi_t)}Vg_x({\xi_t}) , t<T \right)\,.
\end{eqnarray*}
Then, recall Cram\'er's estimate (Lemma \ref{L0}), 
that $\tilde n(T>t)<\infty$ for all $t>0$ and from the first part of Lemma \ref{L3}, that 
$\sup_{x,y}\e^{-\theta(x+y)}|Vg_x(y)|<\infty$.
 Lemma \ref{L3} now entails by dominated convergence that the preceding quantity converges as $x\to -\infty$ towards
$$  c \,
\tilde n\left(F(\xi\circ k_t)  \tilde h({\xi_t}), t<T \right)
\int_0^{\infty} \e^{-\theta z} g(z) \d z =
  c \tilde E\up(F(\xi\circ k_t)) \int_0^{\infty} \e^{-\theta z} g(z) \d z\,,
$$
where $c>0$ is some constant, $\tilde h(y) = \tilde P_{y}(\inf \xi >0)$
and  the identity follows from \eqref{Eq10}. 

Putting the pieces together, we have thus shown that
\begin{eqnarray*}
&&\lim_{x\to -\infty} E_x\left(F(\overleftarrow{\xi}_s, 0\leq s \leq t) g(\sup \xi ), \sigma>t \, \mid \sup Ê\xi  >0\right) \\
&=& c \tilde E\up(F(\xi\circ k_t)) \int_0^{\infty} \e^{-\theta z} g(z) \d z\,.
\end{eqnarray*}
Taking $F\equiv 1$ and invoking Lemma \ref{L0}, we realize that $c=\theta$. 
 This completes the proof.  \QED
 
 \end{section}
 
\begin{section}{Some applications}
We conclude this work by presenting applications of the two main limit theorems,
first to time-reversal identities, then to insurance risk and ruin, and finally to self-similar Markov processes.

\subsection{ Time-reversal identities}
Roughly speaking,  Theorems \ref{T1} and \ref{T2} describe the same two-sided process up-to some time-shift. More precisely, 
consider a two-sided path $\bar\omega$ that has a finite positive maximum and limit $-\infty$ at $\pm \infty$; it should be plain that
the two-step transformation that consists of shifting $\bar \omega$ first at the first instant $\tau$ when it enters $(0,\infty)$ and then shifting again at the first instant when the new path reaches its overall maximum, is equivalent to shifting directly $\bar \omega$ to the first instant $\sigma$ when it reaches its overall maximum. We deduce that
the law under $\p$ of $\xi_{\sigma+\cdot}$,  the path shifted at the instant $\sigma$ when it attains its maximum,  coincides with $\q$, i.e.
\begin{equation}\label{Eq12}
{\mathcal L}(\xi_{\sigma+\cdot}, \p) = \q\,.
\end{equation} 
Conversely, 
the law under $\q$  of $\xi_{\tau+\cdot}$,  the process shifted at the first instant $\tau$ when it enters $(0,\infty)$,  coincides with $\p$, 
\begin{equation}\label{Eq12'}
{\mathcal L}(\xi_{\tau+\cdot}, \q) = \p\,.
\end{equation} 

We now derive a few consequences of these identities. The first concerns  $\overleftarrow{\xi}$, the reversed pre-maximum process.  

\begin{corollary} \label{C3} 
The distribution of $(\overleftarrow{\xi}_t= \sup \xi - \xi_{(\sigma-t)-}, 0\leq t < \sigma)$ under $P$ is the same as that of the process $(\tilde \eta\up_t,0\leq  t<\ell(\epsilon))$, 
where $\tilde \eta\up$ has the law $\tilde P\up$,  $\epsilon$ is an independent exponential variable with parameter $\theta$ and $\ell(\epsilon)=\sup\{t\geq 0: \tilde \eta\up_t \leq \epsilon\}$ denotes  the last passage time of $\tilde \eta\up$  below $\epsilon$.
\end{corollary}

\proof Indeed the law of $(\overleftarrow{\xi}_t)_{0\leq t<\sigma}$ under $P_y$ does not depend on $y$, and as a consequence, it is the same under $P$  as under $\p$. So the claim  is immediate from \eqref{Eq12} and the construction of $\q$.   Alternatively, the interested reader may wish to check directly this identity from Lemma \ref{L3}
 and \eqref{Eq10}.  \QED
 
\noindent {\bf Remark.} {\sl We see from Corollary \ref{C3}  that the distribution of $\sup \xi $ under $P$ is stochastically dominated by the exponential distribution with parameter $\theta$, a property that can also be verified directly from exponential tilting.  It turns out that the law of $\sup \xi $ under $P$ is in fact a divisor of the exponential distribution 
 in the following sense. Let $A$ be a variable distributed as $\sup \xi $ under $P$
 and $B$ an independent variable distributed according to the second marginal of $\rho$,
 where  $\rho$ is the probability law on $\R_+^2$ defined in \eqref{Eq7}. Then $A+B$ has the exponential law with parameter $\theta$. 
Indeed, we see from the construction of $\p$ that the law of the one-sided path
 $(\xi_t)_{t\geq 0}$ under $\p$ coincides with that of $(B+ \xi_t)_{t\geq 0}$ under $P$, assuming implicitly that $B$ and $\xi$ are independent. The claim is now plain from \eqref{Eq12} and the fact that under $\q$, $\sup \xi = \epsilon$ has the exponential distribution with parameter $\theta$. }

We next present a second  application to time-reversal. In this direction, we introduce the law $P'$  of $-\xi$ under $\tilde P$, and observe that $P'$ also fulfills Cram\'er's condition with the same exponent $\theta$, and that the exponentially tilted probability measure associated to $P'$ is the law $\hat P$ of $-\xi$ under $P$. 
We also require the analog of \eqref{Eq3} to hold for $P'$, which amounts to
\begin{equation}\label{Eq3'}
E(|\xi_1|)<\infty\,.
\end{equation}
Last, we denote by   $\p'$  the law on $\bar \Omega$ constructed from $P'$ as in section 3.1, and which corresponds to weak limit of the path shifted at its first passage time above $0$ under $P'_x(\cdot \mid \sup \xi >0)$ as $x\to -\infty$, see Theorem \ref{T1}. We point at the following connexion with $\p$ involving last-passage times.

 \begin{corollary}\label{C4} Assume \eqref{Eq1}, \eqref{Eq2}, \eqref{Eq3} and \eqref{Eq3'}, and 
 introduce
 $$\ell=\sup\{t\in\R: \xi_t>0\}\,,$$
 the last passage time above $0$. Then we have
$$ {\mathcal L}\left( (\xi_{(\ell-t)-})_{t\in\R}, \p\right) = \p'\,.$$

\end{corollary}

\proof Let
$\q'$ denote the law which is constructed from $P'$ as in section 3.2; it is immediately seen from the discussion preceding the statement that $\q'$ coincides with the law of the time-reversed process $(\xi_{-t-})_{t\in\R}$ under $\q$, 
 viz.
\begin{equation}
\label{Eq0}
{\mathcal L}\left( (\xi_{-t-})_{t\in\R}, \q\right) = \q'\,.
\end{equation}
  We can now combine \eqref{Eq12} and \eqref{Eq0} and deduce that
$$
         {\mathcal L}((\xi_{\sigma-t})_{t\in\R},\p) = {\mathcal L}((\xi_{\sigma+t})_{t\in\R},\p').
$$
      Hence the law of $\gamma:=\sup\{t\geq 0 : \xi_{\sigma+t} >0 \}$ under $\p'$
      coincides with that of $\sigma$ under $\p$ and then, since $\ell=\sigma+\gamma$,
      we have
$$
       {\mathcal L}((\xi_{-t-})_{t\in\R},\p) = {\mathcal L}((\xi_{\ell+t})_{t\in\R},\p'),
 $$
      which yields the identity of the statement. \QED
  
\subsection{An application to insurance risk}
We shall now derive an application of the preceding section to insurance risk. 
The reserves of an insurance company are modeled by a process $R=(R_t)_{t\geq 0}$;
in the classical theory \cite{As2},   $R$ is often assumed to be a compound Poisson process
with drift, where the drift coefficient  represents the premium rate and the negative jumps the claims.
More general models based on L\'evy processes are nowadays largely used in the literature, see for instance \cite{KKM}.

Imagine that ruin occurs (i.e. the reserves become negative at some time), but that then the insurance company is not bankrupted and rather can borrow and survive until debt is reimbursed and reserves become
positive again. We are interested in the total time when the insurance company is indebted,
$$D=\int_0^{\infty}{\bf 1}_{\{R_t<0\}} \d t\,.$$
We aim at determining the asymptotic law of $D$ when the initial reserve $R_0$ is large
and $R$ is given by a spectrally negative L\'evy process, conditionally on the event that ruin occurs. 

In this direction, it is convenient to express the reserve process as $R=-\xi$ and work under a probability measure $P_x$ for which $\xi$ is a spectrally positive L\'evy process. In particular the initial reserve is $R_0=-x>0$
and 
$$D=\int_0^{\infty}{\bf 1}_{\{\xi_t>0\}} \d t =\int_0^{\infty}{\bf 1}_{\{\xi_{\tau+t}>0\}} \d t\,.$$
Theorem \ref{T1}  has the following consequence:

\begin{corollary}\label{C5}  Assume \eqref{Eq1}, \eqref{Eq2}, \eqref{Eq3} and \eqref{Eq3'},
and that $\xi$ has no negative jumps $P$-a.s. Then the law of $D$ under $P_x(\cdot \mid \sup \xi >0)$
has a weak limit as $x\to -\infty$; more precisely
$$\lim_{x\to -\infty}P_x(D\in \d t \mid \sup \xi >0)= \theta E(\xi_t^- \exp(\theta \xi_t)) \frac{\d t}{t}\,,\qquad t>0\,,$$
where $\xi_t^-= (-\xi_t)\vee 0$ denotes the negative part of $\xi_t$. 
\end{corollary}

\proof 
It should be clear from Theorem \ref{T1} that the law of $D$ under $P_x(\cdot \mid \sup \xi >0)$
has a weak limit as $x\to -\infty$, and that the limiting law corresponds to  the distribution under $\p$
of 
$$J:=\int_{-\infty}^{\infty}{\bf 1}_{\{\xi_t>0\}} \d t = \int_0^{\ell}{\bf 1}_{\{\xi_t>0\}} \d t
= \int_0^{\ell}{\bf 1}_{\{\xi_{\ell-t}>0\}} \d t\,,$$
 with 
$\ell=\sup\{t\in\R: \xi_t>0\}$ the last passage time of $\xi$ above $0$.
We can then invoke Corollary \ref{C4} (recall also the notation which is used there) to deduce that $J$ has the same law
under $\p$ as under $\p'$.

The point in working with $\p'$ instead of  $\p$ is that, since $\xi$ has no positive jumps $\p'$-a.s., we have
$\p'(\xi_0=0)=1$, and hence the law of $J$ under $\p'$ is the same as under $P'$. According to a celebrated identity due to Sparre Andersen (see for instance Lemma VI.15 in \cite{Be}), under $P'$, $J$ has the same law as $\sigma$, the instant when the path reaches its maximum. All in all, there are the identities
$${\mathcal L}(J,\p)= {\mathcal L}(J,\p')={\mathcal L}(J,P')={\mathcal L}(\sigma, P')\,.$$

The Laplace transform of the latter distribution is well-known for spectrally negative L\'evy processes; see e.g. Theorem VII.4(i) in \cite{Be}. Specifically, if we denote by $\Phi: [0, \infty)\to [\theta,\infty)$ the Laplace exponent of the killed subordinator
$(\tau_x, x\geq 0)$ of the first passage times under $P'$, i.e.
$$E'(\exp(-a\tau_x))= \exp(-x\Phi(a))\,, \qquad a,x\geq 0\,,$$
then
$$E'(\exp(-a\sigma))= \theta/\Phi(a)\,, \qquad a\geq 0\,.$$
Observe that we can re-write 
$$ \frac{\theta}{\Phi(a)}= \theta\int_0^{\infty} \int_0^{\infty} \e^{-at} P'(\tau_y\in \d t) \d y\,.$$
We then use a version of the Ballot Theorem (Corollary VII.3 in \cite{Be}) to express the right-hand side as
$$\theta\int_0^{\infty} \int_0^{\infty} \e^{-at} \frac{y}{t} P'(\xi_t\in \d y) \d t\,.$$
This establishes our claim since
$P'(\xi_t\in \d y) = \tilde P(-\xi_t\in \d y) = \e^{-\theta y} P(-\xi_t\in \d y)$.  \QED

\noindent{\bf Remark.} {\sl It can be checked directly that the limit distribution in Corollary \ref{C3} is a probability measure. 
Indeed, we 
have from the Ballot Theorem that
$$\int_0^{\infty} E(\xi_t^- \exp(\theta \xi_t)) \frac{\d t}{t} =\int_0^{\infty} \tilde E(\xi_t^- ) \frac{\d t}{t} =  \int_0^{\infty}P'(\tau_y<\infty)\d y\,.$$
Since under $P'$ the process $(\tau_y)_{y\geq 0}$ of the first passage times is a subordinator killed at rate $\theta$,   the right-hand side above equals $1/\theta$, showing that the limit distribution  in Corollary \ref{C5} is indeed a probability measure.}
  \subsection{Applications to self-similar Markov processes}
 
 This section is devoted to applications of Theorems \ref{T1} and \ref{T2} to self-similar Markov processes (in short, ssMp); we start by briefly recalling well-known facts due to  Lamperti \cite{La} in this area.
 
 We consider a Markov process $X=(X_t)_{t\geq 0}$ with values in $[0,\infty)$ and denote by
 $\P_y$ its law started from $y>0$. We suppose self-similarity in the sense that 
 $${\mathcal L} ((cX_{t/c})_{t\geq 0}, \P_y) = \P_{cy}\qquad \hbox{for every }c>0 \hbox{ and } y>0\,.$$ 
 Note that for the sake of simplicity, we focus on the case where the index of self-similarity is $1$; of course this induces no loss of generality as we can always replace $X$ by $X^a$ for some adequate power $a\neq 0$. The boundary point $0$ will have a special role; we denote by 
 $$\zeta=\inf\{t\geq 0: X_t=0\}$$ 
 the first hitting time of $0$ and always assume that $X$ reaches $0$ continuously whenever $\zeta < \infty$, that is
 $$\P_y(X_{\zeta-}>0, \zeta<\infty)=0\,.$$
 
 Lamperti  \cite{La} pointed out that ssMp are related to L\'evy processes by a simple transformation that we now describe. The exponential integral
 $$I=\int_0^{\infty}\exp(\xi_t)\d t$$
 of a L\'evy process $(\xi_t)_{t\geq 0}$ converges $P$-a.s. if $\xi$ drifts to $-\infty$, and diverges a.s. otherwise, see e.g. Theorem 1 in \cite{BY}. 
 If we define the time-substitution $\gamma: [0,I)\to \R_+$ by
 $$\int_0^{\gamma(t)}\exp(\xi_s)\d s = t\,,\qquad 0\leq t < I\,,$$ 
 then under $P_x$ the process 
 $$X_t=\exp(\xi_{\gamma(t)}) \quad\hbox{ for }0\leq t < I$$ 
 is a ssMp started from $y=\e^x$
 and (possibly) killed at time $\zeta=I$ when it hits the boundary point $0$. This transformation yields the construction of an arbitrary ssMp on $(0,\infty)$ which reaches $0$ continuously.

 We are interested in the case when $\zeta<\infty$, $\P_y$-a.s., which is equivalent  to assuming that  $\xi$ drifts to $-\infty$ under $P$.  It is then natural to ask whether $X$ possesses a recurrent extension that leaves $0$ continuously. Roughly speaking, we know from the fundamental work of It\=o that this amounts to the existence of an excursion measure ${\bf n}$ under which the process leaves $0$ continuously. 
 This question was raised by Lamperti and solved by him  in the special case where $X$ is a Brownian motion killed at $0$. It was then
 tackled by Vuolle-Apiala \cite{VA} and next by Rivero \cite{Riv1} who obtained sufficient conditions. The problem has been completely solved by Fitzsimmons \cite{Fitz} and 
 Rivero \cite{Riv2} who have shown that the existence of a (unique) recurrent extension that leaves $0$ continuously is equivalent to  Cram\'er's condition \eqref{Eq2} with $\theta<1$ (recall that we assume that the index of self-similarity is $1$).
 
 In the case when furthermore  \eqref{Eq3}  holds, 
 Rivero \cite{Riv1} has obtained some explicit descriptions of the excursion measure ${\bf n}$, including via its entrance law (Theorem 1 in \cite{Riv1}) and {\it \`a la } It\=o, i.e. conditionally on its duration
 (section 4 in \cite{Riv1}). 
 The purpose of this section is to propose further descriptions
 of ${\bf n}$ from the point of view of its height
 $$H=\sup_{0\leq t < \zeta} X_t\,.$$
 We start with the following Lamperti-type representation for the conditional law ${\bf n}(\cdot \mid H>1)$ in terms of $\p$.

\begin{corollary} \label{C1} Assume \eqref{Eq1}, \eqref{Eq2} and \eqref{Eq3}, and that $\theta<1$. The quantity
$$\bar I=\int_{-\infty}^{\infty}\exp( \xi_t) \d t$$
is finite ${\mathcal P}$-a.s. and if we define the time-change $\bar{\gamma}: (0,\bar I)
\to (-\infty,\infty)$ by
$$\int_{-\infty}^{\bar{\gamma}(t)} \exp( \xi_s) \d s = t\,,\qquad 0< t < \bar I\,,$$
then the distribution of the process $(\exp(\xi_{\bar{\gamma} (t)}), 0< t < \bar I)$ under $\p$
coincides with that of $(X_t, 0< t < \zeta)$ under  the conditional law ${\bf n}(\cdot \mid H>1)$. 
\end{corollary}
{\bf Remark.} {\sl The assumption that $\theta<1$ will play no role in the proof; it is only needed in Corollary \ref{C1} to ensure the existence of a true excursion measure. If this assumption is dropped, the statement still holds, except that now ${\bf n}$ is only a pseudo-excursion measure in the sense used by Rivero \cite{Riv1}.} 

\proof  Since $\xi$ drifts to $\infty$ under  $\tilde P$,
the exponential integral $\int_0^{\infty}\exp(-\xi_t)\d t$ is finite $\tilde P_x$-a.s. for every $x\geq 0$, and thus also
almost surely under the conditional laws $\tilde P\up_x$. On the other hand, $\xi$
drifts to $-\infty$ under $P$ and hence $\int_0^{\infty}\exp(\xi_t)\d t <\infty$ $P_y$-a.s. for every $y\geq 0$. By conditioning the bivariate path $(\xi_t)_{t\in\R}$ on $(-\xi_{0-}, \xi_0)$, we now see that $\bar I<\infty$, $\p$-a.s.

Next consider  the ssMp  which is associated by Lamperti's transformation to the tilted L\'evy process with law $\tilde P$. Recall that  \eqref{Eq3} ensures that the expectation  $\tilde E(\xi_1)$ exists in $(0,\infty)$. In turn this guaranties that $0$ is an entrance boundary for the ssMp (see Chaumont {\it et al.} \cite{CKPR} or \cite{BS}), and we denote by  $\tilde \P_0$ the law of the latter starting from $0$.  Then according to the proof of Theorem 1 of Rivero \cite{Riv1} (beware that our notations differ), if $\theta<1$, then the excursion measure ${\bf n}$
 of the recurrent extension of $X$ can be expressed as the Doob's transformation of $\tilde \P_0$ 
 corresponding to the excessive function $x\mapsto x^{-\theta}$, $x>0$. This means that for every stopping time $S$ which is  finite $\tilde \P_0$-a.s., there is the identity
 $${\bf n}(F((X_s)_{ 0\leq s \leq S}), S< \zeta) = \tilde \E_0(F((X_s)_{ 0\leq s \leq S})) X_S^{-\theta})\,,$$
 where $F: \D[0,S]\to \R_+$ is a generic functional that depends on the path only up-to time $S$.
 
 We can now invoke the connexion between the laws $\tilde \P_0$ and $\tilde \p$ via a Lamperti's type transformation described in Corollary 4 of \cite{BS}. In the notation of the present statement,  $\tilde \p(\bar I=\infty)=1$ and 
 the law of $(X_t=\exp(\xi_{\bar{\gamma} (t)}), t > 0)$ under $\tilde \p$ is $\tilde \P_0$.
 It is readily checked in this setting that
 $$S_1:=\inf\{s\geq 0: X_s>1\}= \int_{-\infty}^0 \e^{\xi_s} \d s\,,$$
 and hence $\tilde \p(\bar \gamma(S_1)=0)=1$.
 Applying this construction to any stopping time $S\geq S_1$  which is finite $\tilde \P_0$-a.s. gives 
  $${\bf n}(F((X_s)_{ 0\leq s \leq S}), S< \zeta) = \tilde{\mathcal E} (F((\exp(\xi_{\bar{\gamma} (s)}))_{ 0\leq s \leq S})) \exp(-\theta \xi_{\bar{\gamma}(S)}))\,.$$
By \eqref{Eq6}, the right hand side can be expressed as
 $$\frac{1}{c(\theta)}{\mathcal E} (F((\exp(\xi_{\bar{\gamma} (s)}))_{ 0\leq s \leq S})), \bar{\gamma}(S)<\infty)\,,$$
which establishes our claim, since
$$
         {\bf n}(H>1) = \frac{1}{c(\theta)}\p(\bar\gamma(S_1)<\infty)
                      = \frac{1}{c(\theta)}, 
$$
      where the last equality follows by the local equivalence of $\tilde \p$ and $\p$
      and  $\tilde \p(\bar \gamma(S_1)=0)=1$. \QED

We shall now provide a more precise description of the excursion measure $n$, now conditionally on the height $H$.  In this direction we recall from Corollary 1(ii) of \cite{Riv1} that 
$${\bf n}(H>z)=c z^{-\theta}\,,\qquad z>0\,,$$
where $c>0$ is some constant depending on the choice of the normalization of the excursion measure. Roughly speaking, we shall show that if we  decompose the excursion of $X$ at the instant $\lambda$ when it reaches its maximum, then conditionally on the value of that maximum, the two parts can  be described as Lamperti's type transforms of the underlying L\'evy process conditioned to stay negative and conditioned to stay positive, respectively. More precisely, 
consider two independent processes $\eta\da$ and $\tilde \eta\up$ with respective laws $P\da$ and
$\tilde P\up$. Introduce 
the exponential integrals 
$$I\da=\int_{0}^{\infty}\exp( \eta\da_t) \d t \quad \hbox{and}\quad 
\tilde I\up=\int_{0}^{\infty}\exp(- \tilde \eta\up_t) \d t$$
which are both finite a.s. 
Then define the time-changes $\gamma\da: [0,I\da)
\to [0,\infty)$  and $\tilde \gamma\up : [0,\tilde I\up)
\to [0,\infty)$ by
$$\int_{0}^{\gamma\da(t)} \exp( \eta\da_s) \d s = t 
\qquad \hbox{for}\  
0\leq t < I\da $$
and
$$
\int_{0}^{\tilde \gamma\up(t)} \exp( -\tilde \eta\up_s) \d s = t \qquad \hbox{for}\ 0\leq t < \tilde I\up .$$
Finally define the two-sided path with life-interval $(-\tilde I\up, I\da)$
$$Y_t = \left\{
\begin{matrix}
\exp(-\tilde \eta\up _{\tilde \gamma\up(-t)-}) \, &\hbox{if}& \ -\tilde I\up< t \leq 0 \,,\\
\exp(\eta\da_{\gamma\da(t)})\, &\hbox{if}& \ 0\leq t < I\da\,.
\end{matrix}
\right.$$

\begin{corollary} \label{C2} Assume \eqref{Eq1}, \eqref{Eq2} and \eqref{Eq3}, and that $\theta<1$. In the notation above, for every $y>0$, the distribution of the process
$(yY_{t/y}, -y\tilde I\up< t < y I\da)$
 is a version of that
of $(X_{\lambda+t})_{-\lambda<t<\zeta-\lambda}$  under the conditional law ${\bf n}(\cdot \mid H=y)$, where
$\lambda$ denotes the a.s. unique instant when $X$ reaches its maximum. 
\end{corollary}

\proof  
Since the laws $\p$ and $\q$ are related by a simple time-shift \eqref{Eq12},
 and the time-shift has no effect on the Lamperti's type transformation described in 
 Corollary \ref{C1}, the latter still holds if we replace $\p$ by $\q$. 
 This yields  a Lamperti's type representation for the excursion of the ssMp shifted at the time when it reaches its maximum, conditionally on $H>1$. Specifically, let $\epsilon$ be an exponential variable with parameter $\theta$ which is independent of $\eta\da$ and $\tilde \eta\up$. If we  define
 $$Y^{(\epsilon)}_t = \left\{
\begin{matrix}
\exp(\epsilon -\tilde \eta\up _{\tilde \gamma\up(-t\e^{-\epsilon})}) \, &\hbox{if}& \ - \e^{\epsilon} \tilde I\up< t \leq 0 \,,\\
\exp(\epsilon + \eta\da_{\gamma\da(t\e^{-\epsilon})})\, &\hbox{if}& \ 0\leq t < \e^{\epsilon} I\da\,,
\end{matrix}
\right.$$
then $Y^{(\epsilon)}$ is a version of that
of $(X_{\lambda+t})_{-\lambda<t<\zeta-\lambda}$  under the conditional law ${\bf n}(\cdot \mid H > 1)$, where
$\lambda$ denotes the a.s. unique instant when $X$ reaches its maximum. By conditioning on $\exp(\epsilon)$, this yields our claim when $y>1$, and then 
thanks to the scaling property, for arbitrary $y>0$. \QED

We conclude by noting that this path decomposition {\it \`a la} Williams of the excursion of a ssMp
enlightens a duality identity obtained by Rivero, see Proposition 6(ii) in \cite{Riv1}.
 \end{section}

 \vskip1cm \noindent {\bf Acknowledgment :}  This work has been supported by ANR-08-BLAN-0220-01. It 
  has been undertaken while M. Barczy was on a post-doctoral position at
  the Laboratoire de Probabilit\'es et Mod\`{e}les Al\'eatoires, University Pierre-et-Marie Curie,
  thanks to NKTH-OTKA-EU FP7 (Marie Curie action) co-funded 'MOBILITY' Grant No. MB08-A 81263.
  M. Barczy was also supported by the Hungarian Scientific Research Fund under Grant No.\ OTKA T-079128.

  \end{document}